\documentclass[12pt]{article}
\title{On the structure of homogeneous symplectic varieties of complete intersection}
\author{Yoshinori Namikawa}
\date{ }
\usepackage{amsmath,amssymb,amsthm, amscd}
\chardef\bslash=`\\

\def\0{{\mathcal O}}
\def\g{{\mathfrak g}}
\def\p{{\mathfrak p}}
\def\h{{\mathfrak h}}

\def\k{{\mathfrak k}}

\begin{document}
\maketitle

\begin{center}
{\bf Introduction}
\end{center} 

A normal complex algebraic variety $X$ is called a {\em symplectic variety} (cf. [Be])  
if its regular locus $X_{reg}$ admits a holomorphic symplectic 2-form $\omega$ 
such that it extends to a holomorphic 2-form on a resolution $f: \tilde{X} \to X$. 

Affine symplectic varieties are constructed in various ways such as 
nilpotent orbit closures of a semisimple complex Lie algebra (cf. [CM]),  Slodowy slices to 
nilpotent orbits (cf. [Sl]) or symplectic reductions of  holomorphic symplectic 
manifolds with Hamiltonian actions. Usually these examples come up with 
$\mathbf{C}^*$-actions. 

In this article we shall study a $2n$-dimensional affine symplectic variety $X \subset 
\mathbf{C}^{2n+r}$ defined as a complete intersection of    
$r$ homogeneous polynomials $f_i(z_1, ..., z_{2n+r}) = 0$ $(1 \leq i \leq r)$. Here we assume 
that weights of all coordinates are $1$: $wt(z_1) = ... = wt(z_{2n+r}) =1$. 
The $\mathbf{C}^*$-action on $\mathbf{C}^{2n+r}$ induces a $\mathbf{C}^*$-action 
on $X$. We also assume that the symplectic form $\omega$ is homogeneous 
with respect to this $\mathbf{C}^*$-action. Namely, for some integer $l$, we 
have $t^*\omega = t^l \cdot \omega$ where  $t \in \mathbf{C}^*$. The integer 
$l$ is called the weight of $\omega$ and is denoted by $wt(\omega)$. 
When $X$ is smooth,  $(X, \omega)$ is isomorphic to $(\mathbf{C}^{2n}, \omega_0)$,  
where $\omega_0$ is the standard symplectic 2-form $\Sigma dz_{2i-1} \wedge dz_{2i}$. 
In the remainder we restrict ourselves to the case when $X$ is singular.  

A main result (Main Theorem) is that such an $X$ is isomorphic (as a $\mathbf{C}^*$-variety) 
to the nilpotent variety $N$ of a semisimple complex Lie algebra $\g$ and 
$\omega$ corresponds to the Kostant-Kirillov form.   

The proof consists of two steps. At first we prove that $X$ coincides with a 
nilpotent orbit closure $\bar{O}$ of a semisimple complex Lie algebra $\g$
(Theorem 2). Theorem 2 actually shows that $X$ is the closure of a Richardson orbit $O$ 
and $\bar{O}$ has a crepant resolution. We next prove in {\bf 6} that such a 
nilpotent orbit closure $\bar{O}$ must be the nilpotent variety $N$ if it has complete intersection singularities.    
 
A symplectic variety tends to have a large embedded 
codimension. The main theorem shows that the $A_1$ surface singularity is 
a unique homogeneous symplectic hypersurface.   
As is studied in [LNSV] we have some examples of quasihomogeneous 
symplectic hypersurfaces in higher dimensions. 

The results of this article are concerned with symplectic varieties. However the proof 
of Theorem 2 is based on contact geometry. In particular, a structure theorem 
[KPSW] on contact projective manifolds plays a crucial role.   
We shall remark in the last section that the contact geometry can be also used 
to give another proof of the main theorem of [F].   

The author thanks M. Lehn and C. Sorger for a lot of discussion on 
symplectic hypersurfaces. He also thanks S. Helmke for suggesting to him the approach 
(6.3) by pointing out that the adjoint representation is the lowest dimensional nontrivial irreducible representation in the $E_8$ case.   
\vspace{0.2cm}

{\bf 1}.  Let $X$ be a homogeneous symplectic variety of complete intersection 
defined in Introduction. 

When $X$ is smooth, the polynomials $f_i$ are all linear 
forms; hence we may assume that $r = 0$ and $X = \mathbf{C}^{2n}$. 
We can write  $\omega^n := \omega \wedge ... \wedge \omega = g \cdot 
dz_1 \wedge ... \wedge dz_{2n}$ with a nowhere vanishing homogeneous polynomial 
$g$. Since such a polynomial $g$ must be a constant, we have $wt(\omega) = 2$. 
Now $\omega$ has a form $\Sigma a_{ij}dz_i \wedge dz_j$ with some constants 
$a_{ij}$. Then $\omega$ becomes the standard symplectic 2-form $\Sigma_{1 \le i \le n} 
dz_{2i-1} \wedge dz_{2i}$ after a suitable linear transformation of $\mathbf{C}^{2n}$.    

From now on we consider the case when $X$ is singular. 
Without loss of generality we may assume that $\mathrm{deg}(f_i) \geq 2$ for all $i$. 
In the remainder we put $a_i := \mathrm{deg}(f_i)$. 
By the adjunction formula (or the residue formula) we have 
$$ \omega^n = c \cdot \mathrm{Res}_X(dz_1 \wedge ... \wedge dz_{2n+r}/(f_1, ..., f_r))$$
with a nonzero constant $c$;   
hence   
$$wt(\omega^n) = 2n+r - \Sigma a_i.$$ 

Since $wt(\omega^n) = n\cdot wt(\omega)$ and $wt(\omega) >0$ (cf. [LNSV], Lemma 2.2), 
we have $\Sigma a_i = n+r$ and $wt(\omega) = 1$. 
\vspace{0.2cm}

{\bf Theorem 1}. {\em $X$ has a $\mathbf{C}^*$-equivariant crepant resolution $\pi: Y \to X$.} 
\vspace{0.2cm}

{\em Proof}. 
Let us take a resolution $g: W \to X$ and apply the minimal model program to 
$g$ ([BCHM]). We then finally get  a {\bf Q}-factorial terminalisation 
$\pi: Y \to X$ of $X$. Namely $Y$ has only {\bf Q}-factorial terminal singularities 
and $K_Y = \pi^*K_X$. 
We shall prove that $Y$ is actually smooth. 

The pullback $\pi^*\omega$ defines a symplectic structure on 
the regular part of $Y$. Let $f: Z \to Y$ be a resolution of $Y$. By the assumption 
$(\pi\circ f)^*\omega$ extends 
to a holomorphic 2-form on $Z$; hence $Y$ is a symplectic variety. Then  
$\mathrm{Sing}(Y)$ has even codimension by Kaledin [Ka].  On the other hand,   
since $Y$ has only terminal singularities, $\mathrm{Codim}_Y\mathrm{Sing}(Y) \geq 3$. 
Hence we have $\mathrm{Codim}_{Y}\mathrm{Sing}(Y) \geq 4$. 
Moreover the $\mathbf{C}^*$-action on $X$ extends to a $\mathbf{C}^*$-action on $Y$ (cf. [Na 1, Proposition A.7]). Note here that a symplectic variety has a natural 
Poisson structure and one can consider its Poisson deformation (cf. [Na 2]).  
Take a Poisson deformation $Y_t$ of $Y$. Then the birational map 
$\pi: Y \to X$ also deforms to a birational map $\pi_t: Y_t \to X_t$, where $X_t$ 
is a Poisson deformation of $X$. If we take the Poisson deformation $Y_t$ general 
enough, then $\pi_t$ is an isomorphism (cf. [Na 2, Theorem 5.5]). In particular, $Y_t = X_t$.   
Since $X$ has only complete intersection singularities, so does $Y_t$. On the other hand, $\mathrm{Codim}_{Y_t} \mathrm{Sing}(Y_t) 
\geq 4$. By a result of Beauville [Be, Proposition 1.4], a symplectic singularity is a complete intersection singularity 
only if its singular locus has codimension $\le 3$.  Therefore, $Y_t$ must be smooth. 
Since $Y$ has only {\bf Q}-factorial terminal singularities, any Poisson deformation of 
$Y$ is locally trivial as a flat deformation by Proposition A.9 and Theorem 17 of [Na 1]. 
This means that $Y$ is smooth. Q.E.D.       
\vspace{0.2cm}    

Let us consider the projectivisation $\mathbf{P}(X) := X - \{0\}/\mathbf{C}^*$ of 
$X$. Then the normal projective variety $\mathbf{P}(X)$ admits a contact structure 
with the contact line bundle $O_{\mathbf{P}(X)}(1)$ (cf. [LeB], [Na 3, Section 4]).   
More precisely there is an exact sequence of vector bundles on $\mathbf{P}(X)_{reg}$: 
$$ 0 \to D \to \Theta_{\mathbf{P}(X)_{reg}}  \stackrel{\eta}\to  
O_{\mathbf{P}(X)}(1)\vert_{\mathbf{P}(X)_{reg}} \to 0,$$  
where $\mathrm{rank}(D) = 2n -2$ and $d\eta\vert_D$ induces a non-degenerate 
pairing on $D$. 

{\bf 2}.  We first claim that $\mathbf{P}(X)$ 
also has a crepant resolution \footnote{This is a crucial conclusion obtained from 
the assumption  $wt(z_i)$ are all $1$.}. Let $L$ be a $\pi$-ample line bundle on $Y$. 
If necessary, replacing $L$ by its suitable multiple, we may assume that $L$ 
has a $\mathbf{C}^*$-linearisation (cf. [CG] Theorem 5.1.9). We put $A_m := \Gamma(Y, L^{\otimes m})$ 
for each $m \geq 0$. Note that each $A_m$ has a grading determined by the 
$\mathbf{C}^*$-action. In particular, $A_0$ is the coordinate ring of $X$ and 
$\mathbf{P}(X) = \mathrm{Proj}(A_0)$. Since $A_m$ are graded $A_0$-modules, 
we can consider the associated coherent sheaves $\tilde{A}_m$ on $\mathbf{P}(X)$. 
Define $Z := {\mathcal Proj}_{\mathbf{P}(X)}(\oplus \tilde{A}_m)$.  
Then $Z$ can be identified with $Y - \pi^{-1}(0)/\mathbf{C}^*$ and the 
projective morphism $\bar{\pi}: Z \to \mathbf{P}(X)$ can be identified with 
the natural map $Y - \pi^{-1}(0)/\mathbf{C}^* \to X - \{0\}/\mathbf{C}^*$ 
induced by the $\mathbf{C}^*$-equivariant resolution $\pi: Y \to X$. 
In particular, $\bar{\pi}$ is a birational map.  
Look at the commutative diagram      
\begin{equation} 
\begin{CD} 
Y - \pi^{-1}(0) @>>> Y - \pi^{-1}(0)/\mathbf{C}^* \\ 
@VVV @VVV \\ 
X - \{0\} @>>> X - \{0\}/\mathbf{C}^*.      
\end{CD} 
\end{equation}

Pick a point $x := (z_1(x), ..., z_{2n+r}(x)) \in X - \{0\}$. We have  
$z_i(x) \ne 0$ for some $i$. Define $U_x := X \cap \{(z_1, ..., z_{2n+r}) \in \mathbf{C}^{2n+r}; z_i = z_i(x)\}$. Then 
$U_x$ is isomorphically mapped onto 
a Zariski open subset of $\mathbf{P}(X)$ by the map $X - \{0\} \to 
\mathbf{P}(X)$. 
The map 
$$\sigma_x:  \mathbf{C}^* \times U_x \to X - \{0\}$$ 
sending $(t, x') \in \mathbf{C}^* \times U_x$ to $t\cdot x' \in X -\{0\}$ is an open immersion.  
We put $V_x := \pi^{-1}(U_x)$. Choose a point $y' \in V_x$ 
and put  $x' := \pi (y')$. Denote by 
$O_{x'}$ (resp. $O_{y'}$) the 
$\mathbf{C}^*$-orbit of $x'$ (resp. $y'$). 
  
Since $O_{x'}$ and $O_{y'}$ are both $\mathbf{C}^*$ orbits, there are natural surjections   $\gamma_{x'}: \mathbf{C}^* \to O_{x'}$ $(t \to t\cdot x')$ and 
$\gamma_{y'}: \mathbf{C}^* \to O_{y'}$ $(t \to t\cdot y')$. Moreover 
$\gamma_{y'}$ factorizes $\gamma_{x'}$:  
$$ \mathbf{C}^* \stackrel{\gamma_{y'}}\to O_{y'} \to O_{x'}.$$ 
Since $wt(z_i) = 1$ for all $i$, we see that $\gamma_{x'}$ is an isomorphism; hence 
$\gamma_{y'}$ is also an isomorphism and $O_{y'} \cong O_{x'}$. 

Let  $T_{y'}V_x$ (resp. $T_{y'}O_{y'}$) be the tangent space of 
$V_x$ (resp. $O_{y'}$) at $y'$. Then one has  
$$  T_{y'}V_x \cap T_{y'}O_{y'} = \{0\}.$$ 
In fact, the isomorphism $O_{y'} \to O_{x'}$ induces an isomorphism 
of the tangent spaces $T_{y'}O_{y'} \to T_{x'}O_{x'}$. 
This isomorphism induces an injection $T_{y'}V_x \cap T_{y'}O_{y'} 
\to T_{x'}U_x \cap T_{x'}O_{x'}$. Since $T_{x'}U_x \cap T_{x'}O_{x'} = \{0\}$ 
by the construction of $U_x$, we see that  
$T_{y'}V_x \cap T_{y'}O_{y'} = \{0\}$. 

Let us consider the map 
$$ \sigma_{V_x}: \mathbf{C}^* \times V_x \to Y - \pi^{-1}(0).$$ 
This map induces a map of tangent spaces 
$$T_{(t,y')}(\mathbf{C}^* \times V_x) \to 
T_{t\cdot y'}Y$$ for $(t, y') \in \mathbf{C}^* \times V_x$. 

We claim that $V_x$ is smooth at $y'$ and this map of tangent spaces 
is an isomorphism. We first show the injectivity.    
We identify 
$T_{(t,y')} (\mathbf{C}^* \times \{y'\})$ with $T_t\mathbf{C}^*$ 
and identify $T_{(t,y')}(\{t\} \times V_x)$ with $T_{y'}V_x$. 
Then $T_{(t,y')}(\mathbf{C}^* \times V_x) = 
T_t\mathbf{C}^* \oplus  T_{y'}V_x$. 
Assume that $(\alpha, \beta) \in T_t\mathbf{C}^* \oplus 
T_{y'}V_x$ is sent to zero by the map above. 
The map $\sigma_{V_x}$ induces  isomorphisms  
$\mathbf{C}^* \times \{y'\} \to O_{y'}$  and 
$\{t\} \times V_x \to t\cdot V_x$. Therefore $(\alpha, 0)$ is 
sent to an element of $T_{t\cdot y'}O_{y'}$ and $(0, \beta)$ is sent 
to an element of $T_{t\cdot y'}(t\cdot V_x)$. 
Since $T_{y'}V_x \cap T_{y'}O_{y'} = \{0\}$, we also have 
$T_{t\cdot y'}(t\cdot V_x) \cap  T_{t\cdot y'}O_{y'} = \{0\}$ by the 
$\mathbf{C}^*$-action. This implies that $\alpha = \beta = 0$. 

Note that $\dim Y = \dim V_x +  1$ and $Y$ is smooth. 
If $V_x$ is singular at $y'$, then $\dim T_{y'}V_x > \dim V_x$; 
but then $\dim T_{(t,y')}(\mathbf{C}^* \times V_x)  > 
\dim T_{t\cdot y'}Y$. This contradicts that the above map is an 
injection. Thus $V_x$ must be smooth at $y'$. 
Moreover this implies that the map is an isomorphism. 
  
We finally claim that   
$\sigma_{V_x}$ is an open immersion.   
Assume that two points $(t_i, y_i) \in \mathbf{C}^* \times V_x$, $i = 1, 2$ 
are mapped to the same point of $Y$. Then $y_1$ and $y_2$ are contained in the same 
$\mathbf{C}^*$-orbit. Moreover $\pi (y_1) = \pi (y_2)$. (If $\pi (y_1) \ne \pi (y_2)$, 
then $\pi (y_1)$ and $\pi (y_2)$ must be contained in different $\mathbf{C}^*$-orbits 
because $\sigma_{U_x}$ is an open immersion.)
If $y_1 \ne y_2$, then the natural map $O_{y_1} \to O_{\pi (y_1)}$ of $\mathbf{C}^*$-orbits 
is not a bijection. This contradicts the previous observation. Thus $y_1 = y_2$.  
Then one has $t_1 = t_2$ because 
$\gamma_{y_1}: \mathbf{C}^* \to O_{y_1}$ ($t \to t \cdot y_1$) is an isomorphism. 
This shows that $\sigma_{V_x}$ is an injection. 
Since $\mathbf{C}^* \times V_x$ and $Y$ are both nonsingular and 
the map $T_{(t,y')}(\mathbf{C}^* \times V_x) \to 
T_{t\cdot y'}Y$ is an isomorphism, we see that $\sigma_{V_x}$ is an open 
immersion.  

Now the commutative diagram above is locally identified with 
\begin{equation} 
\begin{CD} 
\mathbf{C}^* \times V_x @>{p_2}>> V_x \\ 
@VVV @VVV \\ 
\mathbf{C}^* \times U_x @>{p_2}>> U_x.      
\end{CD} 
\end{equation} 
By the assumption $\mathbf{C}^* \times V_x \to \mathbf{C}^* \times U_x$ 
is a crepant resolution. This means that $V_x \to U_x$ is also a 
crepant resolution. 
 
Therefore we get a crepant resolution $\bar{\pi}: Z \to \mathbf{P}(X)$ of 
$\mathbf{P}(X)$.  

{\bf  3}. We next claim that $Z$ is a contact projective manifold with the 
contact line bundle $\bar{\pi}^*O_{\mathbf{P}(X)}(1)$. 

For simplicity we write $L$ for $O_{\mathbf{P}(X)}(1)\vert_{\mathbf{P}(X)_{reg}}$. 
The contact structure on $\mathbf{P}(X)_{reg}$ is expressed 
as a twisted 1-form $\eta \in \Gamma (\mathbf{P}(X)_{reg}, 
\Omega^1_{\mathbf{P}(X)_{reg}} \otimes L)$ such that 
$\eta \wedge (d\eta)^{n-1} \in O_{\mathbf{P}(X)_{reg}}$  is nowhere-vanishing. 
In our case $L$ extends to the line bundle $O_{\mathbf{P}(X)}(1)$ on 
$\mathbf{P}(X)$. Let $i: \mathbf{P}(X)_{reg} \to \mathbf{P}(X)$ be the 
natural inclusion map. Since $\mathbf{P}(X)$ has only canonical singularities, 
we have $\bar{\pi}_*\Omega^1_Z \cong i_*\Omega^1_{\mathbf{P}(X)_{reg}}$ 
([GKK]).  
Hence the pull-back $\bar{\pi}^*\eta$ is a section of 
$\Omega^1_Z \otimes \bar{\pi}^*O_{\mathbf{P}(X)}(1)$.  
Moreover, since $\bar{\pi}$ is a crepant resolution, $\bar{\pi}^*\eta 
\wedge (d\bar{\pi}^*\eta)^{n-1}$ is nowhere-vanishing. 

Therefore we get a contact structure of $Z$ with the contact line bundle 
$\bar{\pi}^*O_{\mathbf{P}(X)}(1)$.   

{\bf 4}. When $n = 1$ we already know that $r = 1$ and $f = z_1^2 + z_2^2 + z_3^2$ 
after a suitable change of coordinates (cf. [LNSV], 3.1). 
Note that $Z = \mathbf{P}(X) = \mathbf{P}^1$ 
in this case. We assume that $n \geq 2$. Then $\mathrm{Codim}_X\mathrm{Sing}(X) = 2$ by 
[Be, Proposition 1.4].  
Hence $\mathbf{P}(X)$ actually has singularities and $b_2(Z) \geq 2$. 
Note that $K_Z$ is not nef because $K_Z = \bar{\pi}^*O_X(-n)$. 
By the structure theorem of Kebekus, Peternell, Sommese and Wisniewski [KPSW] 
we conclude that $Z$ is isomorphic to the projectivised cotangent bundle 
$\mathbf{P}(\Theta_M)$\footnote{In this note we employ Grothendieck's notation 
for a projective space bundle. Namely $\mathbf{P}(\Theta_M) = T^*M - (0-section)/\mathbf{C}^*$.} 
of a projective manifold $M$ of dimension $n$; moreover, 
$\bar{\pi}^*O_{\mathbf{P}(X)}(1) \cong O_{\mathbf{P}(\Theta_M)}(1)$. 

Let $\eta_0$ be the canonical contact structure on $\mathbf{P}(\Theta_M)$ 
induced by the canonical symplectic form on $T^*M$. 
Note here that an automorphism $\varphi$ of the vector bundle $\Theta_M$ induces an 
automorphism of $Z:= \mathbf{P}(\Theta_M)$, which is denoted by the same notation 
$\varphi$. Then $\Omega^1_Z$ and $O_{\mathbf{P}(\Theta_M)}(1)$ are both 
$\mathrm{Aut}(\Theta_M)$-linearlized. 
Then our contact form $\eta$ can be written as $\eta = \varphi^*\eta_0$ 
for some $\varphi \in \mathrm{Aut}(\Theta_M)$ (cf. [KPSW], Proposition 2.14).
We may assume that $\eta = \eta_0$
by composing $\varphi$ with the initial identification $Z \cong \mathbf{P}(\Theta_M)$.  

The embedding $X \to \mathbf{C}^{2n+r}$ induces an embedding 
$\mathbf{P}(X) \to \mathbf{P}^{2n+r-1}$. Since $H^0(\mathbf{P}^{2n+r-1}, O_{\mathbf{P}^{2n+r-1}}(1)) 
\cong H^0(\mathbf{P}(X), O_{\mathbf{P}(X)}(1))$, the morphism $\bar{\pi}$ coincides with the 
one defined by the complete linear system $\vert O_{\mathbf{P}(\Theta_M)}(1) \vert$.  
\vspace{0.2cm}

{\bf Lemma}. {\em $\chi (\mathbf{P}(X), O_{\mathbf{P}(X)}) = 1$. }

{\em Proof}.  
We first claim that if $W \subset \mathbf{P}^m$ is a complete intersection 
of type $(d_1, ..., d_k)$, then $\chi (W, O_W(-i)) = 0$  
for all $i  > 0$ with $d_1 + ... + d_k + i < m+1$. We prove this by the induction on $k$. 
Assume that this is true for $k-1$. Let us take the complete intersection $W'$ 
of type $(d_1, ..., d_{k-1})$ such that $W$ is an element of $\vert O_{W'}(d_k) \vert$. 
By the exact sequence 
$$ 0 \to O_{W'}(-i-d_k) \to O_{W'}(-i) \to O_W(-i) \to 0$$ 
we have $\chi (O_W(-i)) = \chi (O_{W'}(-i)) - \chi (O_{W'}(-i-d_k))$. 
Assume that $d_1 + ... + d_k + i < m+1$. Then we have  
$d_1 + ... + d_{k-1} + (i+d_k) < m+1$ and  $d_1 + ... + d_{k-1} + i < m + 1$. 
By the induction assumption $\chi (O_{W'}(-i-d_k)) = \chi (O_{W'}(-i)) = 0$; 
hence $\chi (O_W(-i)) = 0$. 

We next claim that $\chi(W, O_W) = 1$ if $d_1 + ... + d_k  < m+1$. 
This is also proved by the induction on $k$. We take the same $W'$ as above. 
Then $\chi (O_W) = \chi (O_{W'}) - \chi (O_{W'}(-d_k))$. 
By the induction assumption $\chi (O_{W'}) = 1$. By the previous claim 
we have $\chi (O_{W'}(-d_k)) = 0$; hence $\chi (O_W) = 1$ as desired. 

Let us return to the original situation. 
By the argument in {\bf 1} we have $\Sigma a_i < 2n + r$. 
Now one can apply the above claim to $\mathbf{P}(X) \subset \mathbf{P}^{2n+r-1}$. 
Q.E.D.   
\vspace{0.2cm}

Since $\mathbf{P}(X)$ has only rational singularities, we have 
$\chi (Z, O_Z) = \chi (\mathbf{P}(X), O_{\mathbf{P}(X)}) = 1$. 
Let us consider the projection map $p: Z \to M$ of the projective 
space bundle. Since $R^ip_*O_Z = 0$ for $i > 0$, we have 
$\chi (Z, O_Z) = \chi (M, O_M)$. 
In particular, we see that $\chi (M, O_M) = 1$. 

Here we recall a special case of the 
theorem of Demailly, Peternell and Schneider [DPS] 
\vspace{0.2cm}

{\bf Theorem}([DPS, Proposition on p.297]) : 
{\em Let $M$ be a projective manifold with 
nef tangent bundle such that $\chi (M, O_M) \ne 0$.  Then $M$ is a Fano manifold. 
When $\dim M = 2$ or $ 3$, $M$ is a rational 
homogeneous space.} 
\vspace {0.2cm}

In our case we have a much stronger condition. In fact, 
$ O_{\mathbf{P}(\Theta_M)}(1)$ is the pull-back of a very ample line bundle 
by a birational morphism. 
\vspace{0.2cm}

{\bf Proposition}. {\em Let $M$ be a Fano manifold. Assume that 
$\vert O_{\mathbf{P}(\Theta_M)}(1) \vert$ is free from base points.  
Then $M$ is isomorphic to a rational homogeneous 
space, i.e. $M \cong G/P$ with a semisimple complex Lie group $G$ and 
its parabolic subgroup $P$.} 
\vspace{0.2cm}

{\em Proof} .  
The map $H^0(\mathbf{P}(\Theta_M), O_{\mathbf{P}(\Theta_M)}(1)) 
\otimes O_{\mathbf{P}(\Theta_M)} \stackrel{\beta}\to O_{\mathbf{P}(\Theta_M)}(1)$ 
is surjective. 
Let us consider the natural map 
$$ H^0(M, \Theta_M) \otimes O_M \stackrel{\alpha}\to \Theta_M. $$ 
We pull back $\alpha$ by the projection map $p: \mathbf{P}(\Theta_M) 
\to M$. Since $p_*O_{\mathbf{P}(\Theta_M)}(1) = \Theta_M$, 
$p^*\alpha$ factorizes $\beta$: 
$$\beta: H^0(\mathbf{P}(\Theta_M), O_{\mathbf{P}(\Theta_M)}(1)) 
\otimes O_{\mathbf{P}(\Theta_M)} \stackrel{p^*\alpha}\to 
p^*\Theta_M \to O_{\mathbf{P}(\Theta_M)}(1).$$ 
Let $x \in M$ be an arbitrary point and restrict $\beta$ to 
the fibre $p^{-1}(x) \cong \mathbf{P}^{n-1}$. Then we have 
$$\beta (x): H^0(\mathbf{P}(\Theta_M), O_{\mathbf{P}(\Theta_M)}(1)) \otimes 
O_{\mathbf{P}^{n-1}} \stackrel{p^*\alpha(x)}\to 
O_{\mathbf{P}^{n-1}}^{\oplus n} \to O_{\mathbf{P}^{n-1}}(1).$$ 
Note that $\beta (x)$ is also surjective. 
By taking the global sections $\beta (x)$ induces a map 
$\Gamma(\beta (x)): H^0(\mathbf{P}(\Theta_M), O_{\mathbf{P}(\Theta_M)}(1)) \to 
H^0(\mathbf{P}^{n-1}, O_{\mathbf{P}^{n-1}}(1))$. If $\Gamma(\beta (x))$ is not 
surjective, then $\beta (x)$ cannot be surjective. Hence $\Gamma(\beta (x))$ 
must be surjective. This also shows that 
$$\Gamma(p^*\alpha (x)):   
H^0(\mathbf{P}(\Theta_M), O_{\mathbf{P}(\Theta_M)}(1)) \to 
H^0(\mathbf{P}^{n-1}, O_{\mathbf{P}^{n-1}}^{\oplus n})$$ is surjective. 
Since $\Gamma(p^*\alpha (x))$ can be identified with
 the map $H^0(M, \Theta_M) \otimes k(x) \stackrel{\alpha(x)}\to \Theta_M \otimes 
k(x)$, the map $\alpha$ is a surjection by Nakayama's lemma. 

Let $G$ be the neutral component of the automorphism 
group $\mathrm{Aut}(M)$ of $M$. Then $G$ can be written as the extension of 
a complex torus $T$ by a linear algebraic group $L$ (cf. [Fu]) 
$$ 1 \to L \to G \to T \to 1.$$  Note that $q(M) = 0$ because $M$ is a Fano 
manifold. 
If $\dim T > 0$, then $\dim \mathrm{Alb}(M) > 0$  
by Theorem 5.5 of [Fu], which is a contradiction.   
Hence $G$ is a linear algebraic group. As $\alpha$ is surjective, $G$ 
acts transitively on $M$. Therefore $M \cong G/P$ for some parabolic 
subgroup $P$ of $G$ (cf. [Spr, 6.2]). Note that $P$ always contains the radical $r(G)$ of $G$.
Then $r(G)$ acts trivially on $M$; but, since $G$ is the neutral component 
of $\mathrm{Aut}(M)$, $G$ acts effectively on $M$. Hence $r(G) = \{1\}$ 
and $G$ is semisimple.  Q.E.D.  
\vspace{0.2cm}




 

{\bf 5}. 
Assume that $n \geq 2$.    
Now $M$ can be written as $G/P$ with $G$ a semisimple complex Lie group 
and $P$ a parabolic subgroup of $G$. By the proof of the previous proposition 
we may assume that $G = \mathrm{Aut}^0(M)$. 
The cotangent bundle 
$T^*(G/P)$ of $G/P$ has a natural Hamiltonian $G$-action and one can define 
the moment map $\mu: T^*(G/P) \to \g^*$. We identify $\g^*$ with 
$\g$ by the Killing form. Then $\mathrm{Im}(\mu)$ coincides with the closure
$\bar{O}$  of a nilpotent orbit $O \subset \g$. The moment map induces a 
generically finite projective morphism of the projectivisations of $T^*(G/P)$ 
and $\bar{O}$: 
$$ \bar{\mu}: \mathbf{P}(\Theta_{G/P}) \to \mathbf{P}(\bar{O}).$$      
Denote by $O_{\mathbf{P}(\bar{O})}(1)$ the restriction of the tautological line 
bundle $O_{\mathbf{P}(\g)}(1)$ of the projective space $\mathbf{P}(\g)$ to 
$\mathbf{P}(\bar{O})$.  Then it can be checked that 
$O_{\mathbf{P}(\Theta_{G/P})}(1) = \bar{\mu}^*O_{\mathbf{P}(\bar{O})}(1)$
\footnote{Let $\omega_{KK}$ be the Kostant-Kirillov 2-form on 
$O$. Then it gives a contact structure on $\mathbf{P}(O)$ with the contact line bundle 
$O_{\mathbf{P}(O)}(1)$. On the other hand, $\mu^*\omega_{KK}$ is a symplectic 
form on $T^*(G/P)$, which gives a contact structure on $\mathbf{P}(\Theta_{G/P})$ 
with the contact line bundle $\bar{\mu}^*O_{\mathbf{P}(\bar{O})}(1)$. Then we can  
apply [KPSW, Theorem 2.12] to conclude that $\bar{\mu}^*O_{\mathbf{P}(\bar{O})}(1) 
= O_{\mathbf{P}(\Theta_{G/P})}(1)$.}.   
   
This means that $\bar{\pi}: \mathbf{P}(\Theta_{G/P}) \to \mathbf{P}(X)$ 
must be the Stein factorization of $\bar{\mu}$.   

By looking at $\bar{\mu}$ we have an inequality 
$$ (1)\;\; \dim \Gamma (\mathbf{P}(\Theta_{G/P}),   
O_{\mathbf{P}(\Theta_{G/P})}(1)) \geq \dim \Gamma(\mathbf{P}(\bar{O}), 
O_{\mathbf{P}(\bar{O})}(1)).$$

Let $I$ be the ideal sheaf of $\mathbf{P}(\bar{O}) \subset \mathbf{P}(\g)$. 
There is an exact sequence 
$$0 \to H^0(\mathbf{P}(\g), O_{\mathbf{P}(\g)}(1)\otimes I) \to 
H^0(\mathbf{P}(\g), O_{\mathbf{P}(\g)}(1)) \to H^0(\mathbf{P}(\bar{O}), 
O_{\mathbf{P}(\bar{O})}(1)).$$ 
Let $T_0{\bar{O}}$ be the tangent space of $\bar{O}$ at the origin 
$0 \in \bar{O}$. Let $\g = \oplus \g_i$ be the decomposition into the simple 
factors. The closure $\bar{O}$ is the product of nilpotent orbit closures 
$\bar{O_i}$ of $\g_i$. Note that $T_0\bar{O} = \oplus T_0\bar{O}_i$. 
Each $T_0\bar{O}_i$ is a sub $G_i$-representation of the adjoint $G_i$-representation 
of $\g_i$. Since $\g_i$ is an irreducible $G_i$-representation, we have 
$T_0\bar{O}_i = \g_i$. Hence $T_0\bar{O} = \g$.    
This means that there is no hyperplane of $\g$ containing $\bar{O}$; 
hence there is no hyperplane of $\mathbf{P}(\g)$ containing $\mathbf{P}(\bar{O})$.  
This shows that 
$H^0(\mathbf{P}(\g), O_{\mathbf{P}(\g)}(1)\otimes I) = 0$. 
Since $h^0(\mathbf{P}(\g), O_{\mathbf{P}(\g)}(1)) = \dim \g$, we have an inequality 
$$ (2)\;\; \dim \Gamma(\mathbf{P}(\bar{O}), 
O_{\mathbf{P}(\bar{O})}(1)) \geq \dim \g.$$

By (1) and (2) we have an inequality 
$$\dim \Gamma (\mathbf{P}(\Theta_{G/P}),   
O_{\mathbf{P}(\Theta_{G/P})}(1)) \geq \dim \g.$$ 

Since $\Gamma (\mathbf{P}(\Theta_{G/P}),   
O_{\mathbf{P}(\Theta_{G/P})}(1)) = \Gamma (G/P, \Theta_{G/P})$, this inequality 
is actually an equality. 
Hence $\bar{\pi}$ coincides with $\bar{\mu}$  and we have an isomorphism of 
polarised varieties $(\mathbf{P}(X), O_{\mathbf{P}(X)}(1)) \cong (\mathbf{P}(\bar{O}), O_{\mathbf{P}(\bar{O})}(1))$. As $X = \mathrm{Spec}\oplus_{m \geq 0}H^0(\mathbf{P}(X), 
O_{\mathbf{P}(X)}(m))$ and $\bar{O} = \mathrm{Spec}\oplus_{m \geq 0}H^0(\mathbf{P}(\bar{O}), O_{\mathbf{P}(\bar{O})}(m))$, this implies that $X = \bar{O}$. 

 

Finally we give an intrinsic characterization of $G$. 
Notice that we have taken an isomorphism $Z \cong \mathbf{P}(\Theta_M)$ 
such that the contact structure corresponds to the canonical one induced by 
the canonical 2-form on $T^*M$.  
Then $G$ acts on $Z$ as contact automorphisms. Since $\bar{\pi}$ is 
$G$-equivariant, this also means that $G$ acts on $\mathbf{P}(X)_{reg}$ as 
contact automorphisms. The $G$-action determines an embedding  $\g \subset 
H^0(\mathbf{P}(X)_{reg}, \Theta_{\mathbf{P}(X)_{reg}})$.  
  
On the other hand, by [LeB] the contact structure 
$$\Theta_{\mathbf{P}(X)_{reg}}  \stackrel{\eta}\to  
O_{\mathbf{P}(X)}(1)\vert_{\mathbf{P}(X)_{reg}} \to 0$$ 
has a splitting (as $\mathbf{C}$-modules) 
$$s: O_{\mathbf{P}(X)}(1)\vert_{\mathbf{P}(X)_{reg}} \to 
\Theta_{\mathbf{P}(X)_{reg}}$$ so that the subspace 
$$s(H^0(\mathbf{P}(X)_{reg}, O_{\mathbf{P}(X)}(1)\vert_{\mathbf{P}(X)_{reg}})  
\subset H^0(\mathbf{P}(X)_{reg}, \Theta_{\mathbf{P}(X)_{reg}})$$ 
is the infinitesimal contact automorphism group of $\mathbf{P}(X)_{reg}$. 
By the observation above it has the same dimension as $\dim \g$. 
Hence $\g \subset 
H^0(\mathbf{P}(X)_{reg}, \Theta_{\mathbf{P}(X)_{reg}})$ coincides with 
the infinitesimal contact automorphism group of $\mathbf{P}(X)_{reg}$ 
(or $\mathbf{P}(X)$) and $G$ is the neutral component of the 
contact automorphism group 
of $\mathbf{P}(X)$. 

We have thus proved: 
\vspace{0.2cm}

{\bf Theorem 2}. {\em Let $X$ be a singular symplectic variety embedded in 
an affine space $\mathbf{C}^N$ as a complete intersection of homogeneous polynomials. 
Then $X$ coincides with a nilpotent orbit closure $\bar{O}$ of a semisimple complex  
Lie algebra $\g$.}        
\vspace{0.2cm}

By the proof such an orbit $O$ is a Richardson orbit and the Springer 
map $T^*(G/P) \to \bar{O}$ is a birational map.  
 
A typical example of $\bar{O}$ is the nilpotent variety $N$ of $\g$. 
Let $\chi : \g \to \g//G = \mathbf{C}^r$ be the adjoint quotient map. 
Then $N = \chi^{-1}(0)$.   
In particular, $N$ is a complete intersection of $r$ homogeneous polynomials in 
$\g$.   
\vspace{0.2cm} 

The following is the main theorem of this article. 
\vspace{0.2cm}

{\bf Main Theorem}. {\em Let $(X, \omega)$ be a singular symplectic variety embedded in 
an affine space $\mathbf{C}^N$ as a complete intersection of homogeneous polynomials. 
Assume that $\omega$ is also homogeneous. 
Then $(X, \omega)$ coincides with the nilpotent variety $(N, \omega_{KK})$ of a semisimple 
complex Lie algebra $\g$ together with the Kostant-Kirillov form $\omega_{KK}$.} 
\vspace{0.2cm} 

{\bf 6}. 
In this section we prove that the nilpotent orbit closure $\bar{O}$ in Theorem 2 is 
actually the nilpotent variety $N$.  

(6.1) Let $\mathbf{C}[x_1, ..., x_n]$ be a polynomial ring with $n$ variables. For a homogeneous 
ideal $I$ of  $\mathbf{C}[x_1, ..., x_n]$, we put $R := \mathbf{C}[x_1, ..., x_n]/I$ and $d := \dim R$. 
Assume that $I$ does not contain a non-zero homogeneous polynomial of degree $1$. 
We denote by $M$ the maximal ideal $(x_1, ..., x_n)$ of $R$.  

{\bf Lemma}. {\em The following are equivalent.}

(i) {\em The formal completion $\hat{R}$ along $M$ is of complete intersection.} 

(ii) {\em  The ideal $I$ is generated by $n-d$ homogeneous elements.}

{\em Proof}. Since it is clear that (ii) implies (i), we only have to prove that (i) implies  
(ii). The number of minimal generators of $\hat{I}$ equals $\dim_{\mathbf C} (I/IM)$ by 
Nakayama's lemma. The condition (i) then means that $\dim_{\mathbf C}(I/IM) = n-d$. 
One can take $n-d$ homogeneous elements $f_1$, ..., $f_{n-d}$ from $I$ such that 
$\bar{f}_1$, ..., $\bar{f}_{n-d} \in I/IM$ form a basis of $I/IM$. Then it can be 
checked that $f_1$, ..., $f_{n-d}$ actually generate $I$ (cf. the proof of Lemma (A.4) of [Na 1]). 
Q.E.D.     \vspace{0.2cm}

(6.2) Let $R$ be the same as in (6.1) and put $X := \mathrm{Spec}(R)$. 
Assume that a reductive Lie group $G$ acts on $\mathbf{C}^n = \mathrm{Spec}\mathbf{C}[x_1, ..., x_n]$ 
so that $X$ is preserved by $G$. Moreover we assume that 
the $G$-action commutes with the $\mathbf{C}^*$-action on $\mathbf{C}^n$. 

{\bf Lemma}. {\em There are a $G$-representation $V$ with $\dim V = n-d$ and a 
$G$-equivariant morphism $f: \mathbf{C}^n \to V$ of affine spaces such that $f^{-1}(0) = X$.} 

{\em Proof}. Let $I_k$ be the degree $k$ part of the homogeneous ideal $I$. 
Since $G$ respects the grading of $\mathbf{C}[x_1, ..., x_n]$, each $I_k$ is a 
$G$-representation. Let $k_1$ be the minimal number such that $I_{k_1} \ne 0$. 
Let $k_2$ be 
the minimal number $k > k_1$ such that $I'_k := \mathrm{C}[x_1, ..., x_n]_{k-k_1} \cdot I_{k_1}$ 
does not coincide with $I_k$. Since $I'_{k_2}$ is a $G$-subrepresentation 
of $I_{k_2}$, there is a $G$-subrepresentation $I''_{k_2}$ of $I_{k_2}$ such that 
$I_{k_2} = I'_{k_2} \oplus I''_{k_2}$. We next put $I'_k :=     
\mathrm{C}[x_1, ..., x_n]_{k-k_2} \cdot I_{k_2}$ for $k > k_2$ and let $k_3$ be 
the minimal number $k$ such that $I'_k \ne I_k$. Let $I''_{k_3}$ be a $G$-subrepresentation 
of $I_{k_3}$ such that $I_{k_3} = I'_{k_3} \oplus I''_{k_3}$. We repeat this process; then  
$I_{k_1} \oplus I''_{k_2} \oplus I''_{k_3} \oplus ...$ 
becomes a $G$-representation of dimension $n-d$. The $V$ is its dual representation. 
Q.E.D.  \vspace{0.2cm}

(6.3) {\bf Proposition}. {\em A nilpotent orbit closure $\bar{O}$ of an exceptional 
simple Lie algebra $\g$ is of complete intersection if and only if $\bar{O} = N$.} 

{\em Proof}. We put $m:= \dim \g$ and $2n := \dim \bar{O}$. Then $\bar{O}$ is an affine subvariety of 
$\mathbf{C}^m$ with codimension $r := m-2n$. Assume that $\bar{O}$ is defined by $r$ homogeneous 
polynomials $f_i$ with $\mathrm{deg}(f_i) = a_i$. 
As remarked at the beginning of {\bf 1}, we have $\Sigma_{1 \le i \le r}a_i = n + r$. 
Since $a_i \geq 2$ for all $i$, we see that $\Sigma a_i \geq 2r$; thus $n \geq r$. 
In particular, $m = 2n + r \geq 3r$. Therefore we have 
$$\mathrm{Codim}_{\g}\bar{O} \le 1/3 \cdot \dim \g.$$ 
On the other hand, by the previous lemma there are a $G$-representation $V$ with 
$\dim V =  \mathrm{Codim}_{\g}\bar{O}$ and a 
$G$-equivariant map $f: \g \to V$ such that 
$f^{-1}(0) = \bar{O}$. 
There are very few (nontrivial) irreducible representations $V$ of an exceptional simple Lie group 
$G$ with $\dim V < \dim G$ (cf. [F-H], Exercise 24.52 (p.414, see also pp.531,532). 
These are: 

$G_2$: $\dim \g = 14$, $\dim V_{\omega_1} = 7$, 

$F_4$: $\dim \g = 52$, $\dim V_{\omega_4} = 26$,

$E_6$: $\dim \g = 78$, $\dim V_{\omega_1} = \dim V_{\omega_6} = 27$,      
 
$E_7$: $\dim \g = 133$, $\dim V_{\omega_7} = 56$ 

Here we denote by $V_{\omega_i}$ the representations $\Gamma_{\omega_i}$ in [F-H].   
As a consequence, we have no irreducible representation $V$ with $\dim V \le 1/3\cdot \dim \g$. 
Let us look at the $G$-equivariant map $f: \g \to V$. Since there is no irreducible $G$-representation 
of dim $\le 1/3\cdot \dim \g$, the $G$-representation $V$ is a direct sum of trivial representations. 
This means that $\bar{O}$ is the common zeros of some invariant polynomials on $\g$ (with respect to 
the adjoint representation). Notice that the nilpotent variety $N$ of $\g$ is the common zeros of 
{\em all} invariant polynomials on $\g$. Since $\bar{O}$ is contained in $N$, we conclude that 
$\bar{O} = N$. Q.E.D.

(6.4) Let $G$ be a semisimple complex Lie group and let $P$ be a parabolic subgroup 
of $G$. Let $O \subset \g$ be the Richardson orbit for $P$. We assume that 
the closure $\bar{O}$ is normal and the Springer map $T^*(G/P) \to \bar{O}$ is birational. 
One can construct a flat deformation of 
$\bar{O}$ in the following way. Details can be found in [Na 4, Section 2]. 
Let $n(\p)$ (resp. $r(\p)$) be the nilradical (resp. solvable radical) of $\p$. 
Let $\h \subset \p$ be a Cartan subalgebra of $\p$ and define $\k (\p) := 
\h \cap r(\p)$. We then have $r(\p) = \k(\p) \oplus n(\p)$. Notice that 
$\bar{O}$ is the $G$-orbit of $n(\p)$: $$\bar{O} = G \cdot n(\p). $$ 
Then $G \cdot r(\p)$ naturally contains $\bar{O}$. Restricting the adjoint 
quotient map $\chi: \g \to \h/W$ to $G \cdot r(\p)$, we have a map 
$$\chi_{\p}: G \cdot r(\p) \to \h/W.$$ 
Let $\nu: \mathcal{X} \to G \cdot r(\p)$ be the normalization map. 
Then the composition map $\mathcal{X} \to \h/W$ factors through 
$\k(\p)/W'$, where $W' \subset W$ is the stabilizer subgroup of $\k(\p)$ 
as a set: 
$$\chi^n_{\p}: \mathcal{X} \to \k(\p)/W'.$$  
By [Na 4, Proposition 2.6] we have $(\chi^n_{\p})^{-1}(0) = \bar{O}$ \footnote{Notice that 
we assume that $\bar{O}$ is normal.} and 
$\chi^n_{\p}$ gives a flat deformation of $\bar{O}$.  
There is a natural $\mathbf{C}^*$-action on $\mathcal{X}$. If 
$(\chi^n_{\p})^{-1}(0) = \bar{O}$ is of locally complete intersection, 
then all fibres $(\chi^n_{\p})^{-1}(\bar{t})$ are also of locally complete intersection by 
the $\mathbf{C}^*$-action. 

(6.5) A fibre of $\chi_{\p}$ has been already studied in [Sl, 4.3]. 
For $t \in \h$ define $Z_G(t) \subset G$ to be the centralizer of 
$t$ in $G$; namely $$Z_G(t) := \{g \in G; Ad_g(t) = t\}.$$ 
Similarly define $Z_{\g}(t) \subset \g$ to be the centralizer of $t$ in $\g$. Note 
that $Z_{\g}(t)$ is a reductive Lie algebra. 
Then $\p_t := \p \cap Z_{\g}(t)$ is a parabolic subalgebra of $Z_{\g}(t)$. 
Let $O_t \subset Z_{\g}(t)$ be the Richardson orbit for $\p_t$.
Take an element $\bar{\bar{t}}$ from the image of the map $\k(\p) \to \h/W$. Then the fibre 
$\chi_{\p}^{-1}(\bar{\bar{t}})$ can be described as follows. Let 
$\{t_1, ..., t_n\}$ be the inverse image of $\bar{\bar{t}}$ by the map 
$\k(\p) \to \h/W$.  Then one has 
$$\chi_{\p}^{-1}(\bar{\bar{t}}) = \bigcup_{1 \le i \le n} \rho_i(G \times^{Z_G(t_i)} \bar{O}_{t_i}),$$
where $\rho_i: G \times^{Z_G(t_i)} \bar{O}_{t_i} \to G\cdot r(\p)$ is a map defined 
by $\rho_i([g, x]) = Ad_g(x)$.  
As remarked in [Sl, p.56, Remark], $\chi_{\p}^{-1}(\bar{\bar{t}})$ is not necessarily 
irreducible. 
However a fibre of $\chi^n_{\p}$ is always irreducible and normal. 
Consider the Brieskorn-Slodowy diagram ([Na 4, p.728 (2)]):
\begin{equation} 
\begin{CD} 
G \times^P r(\p) @>>> \mathcal{X} \\ 
@VVV @VVV \\ 
\k(\p) @>>> \k(\p)/W'      
\end{CD} 
\end{equation}
Here $G \times^P r(\p)$ gives a simultnaneous resolution of the flat family 
$\mathcal{X} \times_{\k(\p)/W'}\k(\p) \to \k(\p)$. 
Take an element $t$ from $\k(\p)$. The fibre of the map $G \times^P r(\p) \to \k(\p)$ 
over $t$ is $G \times^P (t + n(\p))$. 
Notice that $$G \times^P (t + n(\p)) = G \times^P (P \times^{P_t} n(\p_t)) 
= G \times^{P_t} n(\p_t) = G \times^{Z_G(t)}(Z_G(t) \times^{P_t} n(\p_t)).$$
Let $\bar{t} \in \k(\p)/W'$ be the image of $t$ by the map $\k(\p) \to \k(\p)/W'$. 
Then the map $$G \times^P (t + n(\p)) \to \mathcal{X}_{\bar{t}}$$ coincides with 
the map $$G \times^{Z_G(t)}(Z_G(t) \times^{P_t} n(\p_t)) \to G \times^{Z_G(t)}\tilde{O}_t,$$  
where $\tilde{O}_t$ is the normalization of the orbit closure $\bar{O}_t$.
In particular, one has 
$$(\chi^n_{\p})^{-1}(\bar{t}) = G \times^{Z_G(t)} \tilde{O}_t.$$   
Note that $(\chi^n_{\p})^{-1}(\bar{t})$ is locally the product of 
$G/Z_G(t)$ and $\tilde{O}_t$. If the central fibre $(\chi^n_{\p})^{-1}(0)$ 
is locally of complete intersection, then $\tilde{O}_t$ is locally of 
complete intersection.   

(6.6) Fix a Cartan subalgebra $\h$ of $\g$. Let $\Phi$ be the root system for 
$\g$. Choose a base $\Delta$ of $\Phi$. 
Recall that every parabolic subgroup of $G$ is conjugate to a 
standard parabolic subgroup $P_I$ for a subset $I$ of $\Delta$. 
We denote by $L(P_I)$ the Levi subgroup of $P_I$ containing $H$. 
For example, if $I = \emptyset$, then $P_I$ is a Borel subgroup and 
$L(P_I)$ is nothing but the maximal torus $H$ of $G$.  
In the remainder we assume that $P$ is a standard one $P_I$. 
One has $$\k(\p_I) = \{h \in \h; \alpha (h) = 0, \forall \alpha \in I\}.$$
Define $$\k(\p_I)^{reg} := \{h \in \k(\p_I); \alpha (h) \ne 0, \forall \alpha \in 
\Phi - \Phi_I \}, $$ where $\Phi_I$ is the root subsystem of $\Phi$ generated 
by $I$.  
Choose $\beta \in \Delta - I$ and consider the larger parabolic subgroup 
$P_{I \cup \{\beta \}}$. Then $\k(\p_{I \cup \{\beta \}})$ is naturally 
contained in $\k(\p_I)$. We take an element $t_{\beta}$ from $\k(\p_{I \cup \{\beta \}})^{reg}$. 
Notice that $Z_G(t_{\beta}) = L(P_{I \cup \{\beta \}})$. Moreover $P_I \cap Z_G(t_{\beta})$ is a 
parabolic subgroup of $Z_G(t_{\beta})$, which determines a Richardson orbit $O_{t_{\beta}}$ of 
$Z_{\g}(t_{\beta})$. We then have 
$$ (\chi^n_{\p_I})^{-1}(\bar{t}_{\beta}) \cong G \times^{Z_G(t_{\beta})} \tilde{O}_{t_{\beta}}. $$ 

(6.7) {\em Example}.  Let $P_I$ be the standard parabolic subgroup of $SL(5)$ determined by 
the following marked Dynkin diagram, where the white vertices are simple 
roots belonging to $I$: 

\begin{picture}(300,20)
\put(30,0){\circle*{5}}\put(35,0){\line(1,0){20}}
\put(60,0){\circle{5}}\put(65,0){\line(1,0){20}}
\put(90,0){\circle*{5}}\put(95,0){\line(1,0){20}}
\put(120,0){\circle{5}}
\end{picture} 
\vspace{0.2cm}

We have two black vertices. 
Take the 1-st black vertex as $\beta$. Then the semisimple reduction 
$[Z_{\g}(t_{\beta}), Z_{\g}(t_{\beta})]$ is of type $A_2 \times A_1$. 
Moreover $O_{t_{\beta}}$ is the Richardson orbit of the first $A_2$ for 
the parabolic subalgebra corresponding to  

\begin{picture}(300,20)(0,0) 
\put(30,0){\circle*{5}}\put(35,0){\line(1,0){20}}
\put(60,0){\circle{5}}    
\end{picture} 
\vspace{0.2cm}

Next take the 2-nd black vertex as $\beta$. Then 
$[Z_{\g}(t_{\beta}), Z_{\g}(t_{\beta})]$ is of type $A_3$. 
The orbit $O_{t_{\beta}}$ is the Richardson of $A_3$ for the 
parabolic subalgebra corresponding to 

\begin{picture}(300,20)
\put(30,0){\circle{5}}\put(35,0){\line(1,0){20}}
\put(60,0){\circle*{5}}\put(65,0){\line(1,0){20}}
\put(90,0){\circle{5}}
\end{picture} 
\vspace{0.2cm}

Let $O \subset sl(5)$ be the Richardson orbit for $P_I$. 
Assume that $\bar{O}$ is locally of complete intersection. 
Then $\tilde{O}_t$ is locally of complete intersection for any 
$t \in \k(\p_I)$ by (6.5). As above we take the 1-st black vertex as $\beta$ and consider 
the corresponding $O_{t_{\beta}}$. It is then easily checked that 
$\mathrm{Codim}_{\tilde{O}_{t_{\beta}}}\mathrm{Sing}(\tilde{O}_{t_{\beta}}) = 4$. 
By [Be, Proposition 1.4] $\tilde{O}_{t_{\beta}}$ is not locally of complete intersection. 
This is absurd.  
The second choice of $\beta$ also leads us to a contradiction. 
In this case $\mathrm{Sing}(\tilde{O}_{t_{\beta}})$ has codimension 2 in $\tilde{O}_{t_{\beta}}$ and 
Beauville's proposition cannot be used. Instead we use the previous lemma. 
First notice that every nilpotent orbit 
closure in $sl(m)$ is normal; hence $\tilde{O}_{t_{\beta}} = \bar{O}_{t_{\beta}}$. 
By a direct calculation one has $\dim \bar{O}_{t_{\beta}} = 8$ and $\dim sl(4) = 15$. 
Suppose that $\bar{O}_{t_{\beta}}$ is locally of complete intersection. 
As proved in {\bf 5}, $T_0\bar{O}_{t_{\beta}} = sl(4)$; one can apply Lemma (6.1)  
to the embedding $\bar{O}_{t_{\beta}} \subset sl(4)$. Then $\bar{O}_{t_{\beta}}$ is defined as the 
common zeros of $7$ homogeneous polynomials $f_i$ $(1 \le i \le 7)$. We put 
$a_i := \mathrm{deg}(f_i)$. By the argument at the beginning of {\bf 1} we 
have $a_1 + ... + a_7 = 11$. On the other hand, since $a_i \geq 2$ for all $i$, we 
have $a_1 + ... + a_7 \geq 14$. This is a contradiction. 
\vspace{0.2cm}
   
(6.8)  We are now going to prove that when $\g$ is a classical simple Lie algebra, 
the nilpotent orbit closure $\bar{O}$ in Theorem 2 is actually the nilpotent variety 
$N$. We employ the following strategy. We shall derive a contradiction assuming that 
$\bar{O}$ in Theorem 2 is not the nilpotent variety. 
First we construct a flat deformation of $\bar{O}$: 
$\chi^n_{\p}: \mathcal{X} \to \k(\p)/W'$ as in (6.4). The parabolic subalgebra $\p$ corresponds 
to a marked Dynkin diagram for $\g$. As demonstrated in (6.7), we take a suitable 
simple root $\beta$ and the corresponding element $t_{\beta} \in \k(\p)$ (cf. (6.6)). 
We next consider the fibre $(\chi^n_{\p})^{-1}(t_{\beta})$. Then this fibre is isomorphic 
to $G \times^{Z_G(t_{\beta})}\tilde{O}_{t_{\beta}}$. If $\bar{O}$ is of complete intersection, 
then $\tilde{O}_{t_{\beta}}$ is also of complete intersection. 
But $O_{t_{\beta}}$ is a Richardson orbit in a classical simple Lie algebra which is 
smaller than $\g$. Moreover the corresponding parabolic subalgebra (= the polarization 
of $O_{t_{\beta}}$) is a maximal parabolic subalgebra. 
Finally we derive a contradiction in such a case.    
\vspace{0.2cm}
     
We first treat the case $\g$ is of type A.   
\vspace{0.2cm}

{\bf Proposition}. {\em A nilpotent orbit closure $\bar{O}$ of $sl(m)$ has complete intersection 
singularities if and only if $\bar{O} = N$.} 

{\em Proof}. Note that every nilpotent orbit $O$ of $\g := sl(m)$ is a Richardson orbit and 
its closure is normal. Moreover the Springer map $T^*(G/P) \to \bar{O}$ is birational. 
As remarked just above, we only have to prove that $\bar{O}$ does not have 
complete intersection singularities when $P$ is a {\em maximal} parabolic subgroup of 
$SL(m)$ with $m \geq 3$. Namely $P$ corresponds to  
to a marked Dynkin diagram with only one black vertex:  

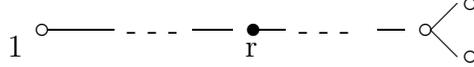
\begin{picture}(300,20)(0,0) 
\put(20, -10){1}\put(30,-3){$\circ$}\put(35,0){\line(1,0){25}} 
\put(65,-3.5){- - -}\put(90,0){\line(1,0){15}}
\put(105,-3){$\bullet$}\put(110,0){\line(1,0){10}}
\put(100,-10){r}\put(125,-3.5){- - -}\put(150,0)
{\line(1,0){55}}\put(207,-3){$\circ$}    
\end{picture}  
\vspace{0.2cm}
   
When $r \ne m/2$, one has $\mathrm{Codim}_{\bar{O}} \mathrm{Sing}(\bar{O}) \geq 4$. 
Then $\bar{O}$ does not have complete intersection singularities by [Be, Proposition 1.4]. 
Assume that $\bar{O}$ has complete intersection singularities when $r = m/2$. 
By a direct calculation we have $\dim \bar{O} = 2r^2$ and $\dim sl(m) = 4r^2 -1$. 
By Lemma (6.1) $\bar{O}$ is a subvariety of $\mathbf{C}^{4r^2 -1}$ defined as 
the complete intersection of $2r^2 -1$ homogeneous polynomials $f_i$. 
We put $a_i := \mathrm{deg}(f_i)$. As discussed at the beginning of {\bf 1}, 
$\Sigma a_i = r^2 + (2r^2 - 1)$. On the other hand, since $a_i \geq 2$, we have 
$\Sigma a_i \geq 2(2r^2 - 1)$. Combining these inequalities we get 
$$ r^2 \geq 2r^2 - 1,$$ which implies that $r = 1$ and then $m = 2$. This contradicts 
the first assumption that $m \geq 3$. Q.E.D.   
\vspace{0.2cm}

(6.9) Let $G$ be $Sp(2n)$ or $SO(n)$ and let $P_I$ be a maximal parabolic subgroup. 
Namely $P$ is the standard parabolic 
subgroup corresponding to one of the following Dynkin diagram.   

$C_n$ 

\begin{picture}(300,20)(0,0) 
\put(30,-3){$\circ$}\put(20,-10){1} \put(35,0){\line(1,0){25}} 
\put(65,-3.5){- - -}\put(90,0){\line(1,0){15}}
\put(105,-3){$\bullet$}\put(110,0){\line(1,0){10}}
\put(100,-10){r}\put(125,-3.5){- - -}\put(150,0)
{\line(1,0){15}}\put(170,-3){$\circ$}\put(175,-3){$\Leftarrow$}
\put(187,-3){$\circ$}    
\end{picture} 
\vspace{0.4cm}

$B_{[n/2]}$ 

\begin{picture}(300,20)(0,0) 
\put(30,-3){$\circ$}\put(20,-10){1} \put(35,0){\line(1,0){25}} 
\put(65,-3.5){- - -}\put(90,0){\line(1,0){15}}
\put(105,-3){$\bullet$}\put(110,0){\line(1,0){10}}
\put(100,-10){r}\put(125,-3.5){- - -}\put(150,0)
{\line(1,0){15}}\put(170,-3){$\circ$}\put(175,-3){$\Rightarrow$}
\put(187,-3){$\circ$}    
\end{picture}
\vspace{0.4cm}

$D_{n/2}$   

\begin{picture}(300,20)(0,0) 
\put(20,-10){1} \put(30,-3){$\circ$}\put(35,0){\line(1,0){25}}
\put(65,-3.5){- - -}\put(90,0){\line(1,0){15}}\put(110,-3)
{$\bullet$}\put(110,-10){r}\put(115,0){\line(1,0){10}}
\put(130,-3.5){- - -} \put(160,0){\line(1,0){10}}
\put(175,-3){$\circ$}
\put(180,0){\line(1,-1){10}}\put(180,0){\line(1,1){10}}
\put(192, 7){$\circ$}\put(192,-13){$\circ$}
\end{picture}
\vspace{0.4cm} 
 
Let $O \subset \g$ be the Richardson orbit for $P_I$. 
We shall prove that $\bar{O}$ is not a homogeneous symplectic variety of complete intersection. 
When $G = Sp(2n)$, the parabolic subgroup $P_I$ is the stabilizer group 
of an isotropic flag of type $(r, 2n-r, r)$. Let $Gr_{iso}(r, 2n)$ be the isotropic Grassmann 
variety parametrizing  such flags. Then  
$$\dim Gr_{iso}(r, 2n)  = \dim Gr(r, 2n) - 1/2\cdot r(r-1) = r(2n-r) - 1/2\cdot r(r-1).$$
Since $\dim \bar{O} = 2\dim Gr_{iso}(r, 2n)$, we have 
$\dim \bar{O} = 2r(2n-r) - r(r-1)$. On the other hand, $\dim sp(2n) = 2n^2 + n$, hence 
$\mathrm{Codim}_{sp(2n)}\bar{O} = 2n^2 + n - 4rn + 3r^2 - r$. 
Assume that $\bar{O}$ is of complete intersection in $sp(2n)$. Let $f_i$ be the defining 
equations of $\bar{O}$ and put $a_i := \mathrm{deg}(f_i)$. Then 
$\Sigma a_i = 1/2\cdot \dim \bar{O} + \mathrm{Codim}_{sp(2n)}\bar{O}$ by 
{\bf 1}. Since $a_i \geq 2$ for all $i$, we have $(3r-2n-1)(3r-2n) \le 0$. 
The only possibilities are following two cases:
 
(i) $n = 3k$ for some integer $k$ and $r = 2k$. 

(ii) $n = 3k + 1$ for some integer $k$ and $r = 2k+1$. 

In both cases $a_i = 2$ for all $i$ (i.e. $\dim V = 1/3 \cdot \dim sp(2n)$.)  
In the first case $O = O_{[3^{2k}]}$ (i.e the nilpotent orbit consisting of 
the matrices of Jordan type $(3, ..., 3)$ ($2k$ Jordan blocks of size $3$).  
In the second case $O = O_{[3^{2k},2]}$.
Assume that $O_{[3^{2k}]} \subset sp(6k)$ is of complete intersection.   
By the calculation above we have $\mathrm{codim}_{sp(6k)}\bar{O} = 6k^2 + k$. 
By Lemma (6.2) there are a $G$-representation $V$ of dim $6k^2 + k$ and a 
$G$-equivariant map $f: sp(6k) \to V$ such that $f^{-1}(0) = \bar{O}$. 
By the construction of $V$ (cf. Lemma (6.2)), the dual representation $V^*$ coincides  
with $I_2$ because $a_i = 2$ for all $i$.  
But there is only one (adjoint) invariant quadratic polynomial on $sp(6k)$ up to 
constant. Hence $V$ contains one and only one trivial representation as a direct factor.  
Since an irreducible representation of $sp(6k)$ with $\dim \le 1/3\cdot \dim sp(6k)$ 
is a trivial representation or a standard representation (cf. [F-H], p.531, (24.52)), $V$ is a direct sum of a 
trivial representation and a finite number of standard representations. 

Let us consider the first case (i).  
Notice that, in this case, $\dim V = 1 + (6k^2 + k -1)$. 
If $k \geq 2$, then $6k$ does not divide $6k^2 + k -1$, which is a contradiction.   
When $k = 1$, one has $\dim V = 7$ and $V$ may possibly be 
a direct 
sum of the 6-dimensional standard representation and the trivial representation. 
Since $a_i = 2$ for all $i$, these irreducible factors must be contained in $\mathrm{Sym}^2(sp(6k)^*)$ 
the 2-nd symmetric product of the dual representation of the adjoint one. By the Killing form 
$\mathrm{Sym}^2(sp(6)^*) \cong \mathrm{Sym}^2(sp(6))$ as $Sp(6)$-representations. 
It is easily checked that $\mathrm{Sym}^2(sp(6))$ does not contain the standard representation as a 
direct factor. Hence we have a contradiction also in this case.    

In the second case (ii) we have $\dim V = 1 + (6k^2 + 5k)$. 
Noticing that the standard representation has dimension $6k + 2$, we 
write $6k^2 + 5k = k(6k + 2) + 3k$; hence $6k + 2$ does not divide 
$6k^2 + 5k$. This is a contradiction. 
   
Assume that $G = SO(n)$ and $\bar{O}$  
has complete intersection singularities.  Since $a_i \geq 2$ for all $i$, 
the equality  
$$\Sigma a_i = 1/2\cdot \dim \bar{O} + \mathrm{Codim}_{so(n)}\bar{O}$$ 
implies that $(3r - n)(3r - n+1) \le 0$. There are two possibilities:  

(i) $n = 3k$ for some integer $k$, $r = k$ and $O = O_{[3^k]}$.　

(ii) $n = 3k + 1$ for some integer $k$, $r = k$ and $O = O_{[3^k, 1]}$. 

In both cases $a_i = 2$ for all $i$ (i.e. $\dim V = 1/3 \cdot \dim so(n)$).  
We can again use Lemma (6.2) to have a $G$-equivariant map $f: so(n) \to V$. 
Put $\g = so(n)$ with $n = 3k$ or $n = 3k + 1$. 
Then $\dim V $ is respectively $1/2\cdot(3 k^2 - k)$ or $1/2\cdot(3k^2 + k)$.   
Note that an irreducible representation of $\g$ with dim $\leq   
1/3 \cdot \dim \g$ is a trivial representation or a standard representation 
(cf. [F-H],  p.531, (24.52): Note that, when $\g$ is of $D_4$, two more 
different irreducible representations exist, but the $D_4$ case is not contained in 
the case (i) or the case (ii).).   
Since there is only one (adjoint) invariant quadratic polynomial on $so(n)$ up to 
constant, $V$ is a direct sum of a 
trivial representation and a finite number of standard representations. 
By writing $k = 2l$ or $k = 2l + 1$ according as $k$ is even or odd,  
one can easily check that $\dim V - 1$ is not divided by $n$ in both cases; 
hence we have a contradiction. 
\vspace{0.2cm}

(6.10) Let $\g$ be a complex simple Lie algebra of type $B$, $C$ or $D$. 
Let $O$ be the Richardson orbit of $\g$ for a parabolic subgroup $P$ of $G$. 
Assume that the Springer map $s: T^*(G/P) \to \bar{O}$ is birational. 

{\bf Proposition} {\em The closure $\bar{O}$ of such an orbit is of complete intersection 
if and only if $\bar{O} = N$.}

{\em Proof}. We only have to deal with a Richadson orbit for a standard 
parabolic subgroup $P_I$.  If the Dynkin diagram corresponding to 
$P_I$ has only one black vertex, then we have already checked that 
$\bar{O}$ is not of complete intersection. 
Assume that there are more than one black vertices, but at least one 
vertex is a white vertex. 
Take a white vertex $w$ on the leftmost position. Note that if the 
Dynkin diagram is of type $B$ or $C$, it is unique, but if the 
Dynkin diagram is of type $D$, the choice of such a vertex might have 
two possibilities. 

If there is a black vertex $b$ left adjacent to $w$, then 
take the simple root $\beta$ corresponding to $b$ and apply (6.6). 
Then the problem is reduced to the case where the Dynkin diagram is 
of type $A$ and has only one black vertex with $r = 1$,  or the Dynkin 
diagram is  a smaller one of the same type as $\g$ and has only 
one black vertex with $r = 1$. In each case $\bar{O}_{t_{\beta}}$ is 
normal; we only have to check this in the second case. There is a 
nilpotent orbit ${O}'_{t_{\beta}} \subset \bar{O}_{t_{\beta}}$ such that 
$\mathrm{Codim}_{\bar{O}_{t_{\beta}}}\bar{O}'_{t_{\beta}} = 2$. 
One can check that $\mathrm{Sing}(\bar{O}_{t_{\beta}}, {O}'_{t_{\beta}})$ 
is of type $a$ or of type $g$ in the list of [K-P, p.551]. By Theorem 1, (b) of 
[K-P] we see that $\bar{O}_{t_{\beta}}$ is normal.  
Moreover, in each case, $\bar{O}_{t_{\beta}}$ is not of   
complete intersection (cf. (6.8), (6.9)). 
By the argument in (6.5), the original nilpotent orbit closure 
$\bar{O}$ is not of complete intersection.   

Assume that  there is no black vertex left adjacent to $w$. 
By the definition of $w$ this means that $w$ is on the leftmost position 
on the diagram. In this case we consider the maximal connected Dynkin subdiagram 
$\mathcal{D}$ containing $w$ whose vertices are all white.  
Let $w'$ be a vertex on the rightest position of $\mathcal{D}$. 
Let $b$ be a black vertex right adjacent to $w'$. We take the simple 
root $\beta$ corresponding to $b$ and apply (6.6). 
Then the problem is reduced to the case where the Dynkin diagram  
is of type $A$ and has only one black vertex. 
Then  $\bar{O}_{t_{\beta}}$ is 
normal and is not of complete intersection (cf. (6.8)). 
By the argument in (6.5), the original nilpotent orbit closure 
$\bar{O}$ is not of complete intersection.  Q.E.D. 

(6.11)  Let $O$ be a Richardson orbit of a complex semisimple Lie algebra 
$\g$. Let $\g = \oplus_{1 \le i \le m} \g_i$ be the decomposition into the simple factors. 
Then we have $\bar{O} =  \bar{O}_1 \times ... \times \bar{O}_m$ where each $O_i$ is a 
Richadson orbit of $\g_i$. If the Springer map $T^*(G/P) \to \bar{O}$ is birational, then 
each Springer map $T^*(G_i/P_i) \to \bar{O}_i$ is birational. 
Assume that $\bar{O}$ is of complete intersection. Then each $\bar{O}_i$ is 
also of complete intersection. By (6.3), (6.8) and (6.10) each $\bar{O}_i$ coincides 
with the nilpotent variety $N_i$ of $\g_i$. Then $\bar{O}$ is the nilpotent variety 
$N$ of $\g$.   
\vspace{0.2cm}    
   

 


{\bf 7}.  
{\bf Remarks}

(1) What happens in Main theorem if we do not assume $\omega$ is homogeneous ?  
The author does not know the answer, but  
the following example would be instructive. Let $X \subset \mathbf{C}^5$ be a 
hypersurface defined by $z_1^2 + z_2^2 + z_3^2 = 0$, where $(z_1, ..., z_5)$ are 
coordinates of $\mathbf{C}^5$. Note that $X = S \times \mathbf{C}^2$, where 
$S \subset \mathbf{C}^3$ is a hypersurface defined by $f: = z_1^2 + z_2^2 + z_3^2 = 0$. 
We put $\omega_S := \mathrm{Res}(dz_1 \wedge dz_2 \wedge dz_3/f)$ and 
$\omega_{\mathbf{C}^2} := dz_4 \wedge dz_5$. Define $\omega :=  
\omega_S + \omega_{\mathbf{C}^2}$. Then $(X,\omega)$ is an affine symplectic 
variety. But $\omega$ is not homogeneous because 
$wt(\omega_S) = 1$ and $wt(\omega_{\mathbf{C}^2}) = 2$. 
Note that $\omega \wedge \omega$ is a holomorphic volume form on $X$ of weight $3$. 
Let us prove that there is no homogeneous symplectic 2-form on $X$. 
Assume that such a form $\Omega$ exists. Then $\Omega \wedge \Omega$ is a 
holomorphic volume form on $X$ of an even weight, say $2m$. 
Then one can write $\Omega \wedge \Omega = g\cdot \omega \wedge \omega$ 
with a nowhere vanishing function $g$ of {\em nonzero} weight. 
But such $g$ does not exists; hence one gets a contradiction.       
\vspace{0.2cm}

(2) Let $X$ be an affine symplectic variety in 
$\mathbf{C}^{N}$ defined by a homogeneous ideal $I$ 
(not necessarily of complete intersection) where 
$I$ contains no nonzero linear form.  Denote by 
$R$ the coordinate ring of $X$. By the assumption $R$ is 
graded: $R = \oplus_{n \geq 0}R_n$.    
Assume that $wt(\omega) = 1$. Then $\omega$ induces a 
Poisson structure on $R$ of weight $-1$. In particular, it induces 
a Lie algebra structure on $R_1$ 
$$[\cdot, \cdot]: R_1 \times R_1 \to R_1.$$ 
Let us call this Lie algebra $\g$. Since $R_1 = T^*_0X$, 
we have $\dim \g = N$. 
The natural surjection
$\oplus\mathrm{Sym}^i(R_1) \to R$ induces a closed embedding 
$X \to \g^*$. To prove that $\g$ is semisimple, it seems that one needs  
some geometric arguments as in {\bf 1} - {\bf 5}. 
When $\g$ is semisimple, $\g^*$ is identified with 
$\g$ by the Killing form. This is nothing but the closed embedding 
$X \to \g$ of Main theorem, where $X$ is identified with the nilpotent 
variety $N$.   
\vspace{0.2cm}
 
(3) Let $X$ be the same as in (2).  
Then $\mathbf{P}(X)$ admits a contact 
structure with the contact line bundle $O_{\mathbf{P}(X)}(1)$ in the 
sense of {\bf 1}.  
Let $G$ be the contact automorphism group of $\mathbf{P}(X)_{reg}$. 
The Lie algebra $\g$ is contained in $H^0(\mathbf{P}(X), \Theta_{\mathbf{P}(X)})$ 
and the map $H^0(\mathbf{P}(X), \Theta_{\mathbf{P}(X)}) \stackrel{\eta}\to 
H^0(\mathbf{P}(X), O_{\mathbf{P}(X)}(1))$ induces an isomorphism 
$\g \cong H^0(\mathbf{P}(X), O_{\mathbf{P}(X)}(1))$ by [Be 2], Proposition 1.1. 
In general we only know that $\dim \g \geq N$.  
The closed embedding $\mathbf{P}(X) \to \mathbf{P}(\g^*)$ is a $G$-equivariant map. 
By a similar argument to [Be 2], Section 1, the $G$-action on $\mathbf{P}(X)$ lifts 
to a $G$-action on $X$. Moreover the above embedding lifts to a 
$G$-equivariant closed embedding $X \to \g^*$. By this embedding $X$ is identified 
with a coadjoint orbit closure of $\g^*$. In particular, $G$ acts transitively  
on $X_{reg}$. But we do not know when $G$ is semisimple.   
\vspace{0.2cm}
 
(4) One can give another proof of [F, Main theorem] : 
  
{\em Every crepant resolution of a 
nilpotent orbit closure $\bar{O}$ of a semisimple complex Lie algebra $\g$ 
coincides with a Springer resolution $\mu: T^*(G/P) \to \bar{O}$.} 

The following proof can be regarded as a translation of the original proof 
into contact geometry. 
 
If $\bar{O}$ has a crepant resolution $\pi: Y \to \bar{O}$, then $\mathbf{P}(\bar{O})$ 
also has a crepant resolution $\bar{\pi}: Z \to \mathbf{P}(\bar{O})$. 
The Kostant-Kirillov 2-form on $O$ induces a contact structure on $\mathbf{P}(O)$.  
Moreover this contact structure is pulled back to a contact structure $\eta$ on $Z$. 
The contact structure $\eta$ can be regarded as an element of  
$\Gamma (Z, \Omega^1_Z \otimes {\bar{\pi}}^*O_{\mathbf{P}(\bar{O})}(1))$. 
Here $O_{\mathbf{P}(\bar{O})}(1)$ 
is the pull-back of $O_{\mathbf{P}(\g)}(1)$ by the inclusion map 
$\mathbf{P}(\bar{O}) \to \mathbf{P}(\g)$. 

Assume that $b_2(Z) > 1$. Then $Z$ is isomorphic to $\mathbf{P}(\Theta_M)$ 
for a projective manifold $M$. 
Let $\eta_0$ be the canonical contact structure on $\mathbf{P}(\Theta_M)$ 
induced by the canonical symplectic form on $T^*M$. 
By the same argument as in {\bf 4}, we may assume that $\eta = \eta_0$. 
Let $G$ be the adjoint group of $\g$. We 
prove that $M \cong G/P$ with some parabolic subgroup $P$ of $G$. 
By [F, Proposition 3.1] the $G$-action on $\bar{O}$ extends to a $G$-action 
on $Y$. Since the $G$-action is compatible with the $\mathbf{C}^*$-action, 
we have a $G$-action on $Z$. 

We first claim that this $G$-action is induced by a $G$-action on 
$M$.  
It is well known that all contact automorphisms of $(\mathbf{P}(\Theta_M), \eta_0)$ 
are those induced by the automorphisms of $M$.  
Since our $G$ acts on $\mathbf{P}(\Theta_M)$ as contact automorphisms,    
our claim has been justified. 

We next claim that the $G$-action on $M$ is transitive. 
Since $O_{\mathbf{P}(\Theta_M)}(-1) = {\bar{\pi}}^*O_{\mathbf{P}(\bar{O})}(-1)$, 
the map $\bar{\pi}$ pulls back the $\mathbf{C}^*$-bundle 
$\bar{O} - \{0\} \to \mathbf{P}(\bar{O})$ to the $\mathbf{C}^*$-bundle 
$T^*M - (0-section) \to \mathbf{P}(\Theta_M)$: 
\begin{equation} 
\begin{CD} 
T^*M - (0-section) @>>> \mathbf{P}(\Theta_M) \\ 
@VVV @VVV \\ 
\bar{O}-\{0\} @>>> \mathbf{P}(\bar{O})      
\end{CD} 
\end{equation}

The $G$-action on $M$ induces a natural $G$-action on 
$T^*M - (0-section)$. It induces a $G$-linearization of $O_{\mathbf{P}(\Theta_M)}(-1)$. 
On the other hand, $\bar{O} - \{0\}$ has 
a natural $G$-action and it induces a $G$-linearization of 
$O_{\mathbf{P}(\bar{O})}(-1)$. The crepant resolution $\bar{\pi}$ is an 
isomorphism over $\mathbf{P}(O)$. By the identification of $\bar{\pi}^{-1}(\mathbf{P}(O))$ 
with $\mathbf{P}(O)$ two line bundles 
$O_{\mathbf{P}(\Theta_M)}(-1)\vert_{\bar{\pi}^{-1}(\mathbf{P}(O))}$ and 
$O_{\mathbf{P}(\bar{O})}(-1)\vert_{\mathbf{P}(O)}$ are identified. 
Each one has a $G$-linearization coming from that of $O_{\mathbf{P}(\Theta_M)}(-1)$ 
or $O_{\mathbf{P}(\bar{O})}(-1)$. By the uniqueness of the $G$-linearization 
([Mu], Proposition 1.4) these two $G$-linearizations are the same. In particular, 
the  commutative diagram above is 
$G$-equivariant.  Since $G$ has an open dense orbit $O$ in $\bar{O}$, 
it also has an open dense orbit in $T^*M - (0-section)$. 
This also shows that $G$ has an open dense orbit $U$ in $M$.
The following argument is the same as in [F]. Write $U = G/P$ with a 
closed subgroup $P$ of $G$. Note that $T^*U = G \times^P (\g/\p)^*$ and 
$G$ has an open dense orbit in $T^*U$. This implies that $P$ has an 
open dense orbit in $(\g/\p)^*$ by the coadjoint action. Then, by Proposition 
3.10 of [F] we see that $P$ is a parabolic subgroup, which implies that 
$U$ is a projective manifold. Since $U$ is an open dense subset of $M$, we 
must have $M = U$. 

Let $\mu: T^*(G/P) \to \g^*$ be the moment map 
and let $\bar{\mu}: \mathbf{P}(\Theta_{G/P}) \to \mathbf{P}(\g^*)$ 
be its projectivization. Note that $\mathrm{Im}(\mu) = \bar{O}' $ with 
a coadjoint orbit $O'$ of $\g^*$. Then $\mathrm{Im}(\bar{\mu}) = 
\mathbf{P}(\bar{O}')$. To compare the map $\bar{\mu}$ with $\bar{\pi}$, 
we identify the nilpotent orbit $O$ with a coadjoint orbit of $\g^*$ by 
the Killing form $\g \cong \g^*$. 

We start with a rather general setting: let $Z$ be a projective contact 
manifold with the contact structure 
$$ 0 \to D \to \Theta_Z \stackrel{\eta}\to L \to 0.$$ 
Assume that a semisimple complex Lie group $G$ acts effectively on $Z$ as 
contact automorphisms. Let $\g \subset H^0(Z, \Theta_Z)$ be the 
space of infinitesimal contact automorphisms determined by 
$G$ and let $V \subset H^0(Z, L)$ be the image of $\g$ by the map 
$H^0(Z, \Theta_Z) \stackrel{\eta}\to H^0(Z, L)$. 
\vspace{0.2cm}

{\bf Lemma}. {\em Let $O \subset \g^*$ be a coadjoint orbit preserved 
by the natural $\mathbf{C}^*$-action of $\g^*$. Assume that 
$$f: Z \to \mathbf{P}(\bar{O})$$ is a generically finite surjective $G$-equivariant 
morphism such that $L = f^*O_{\mathbf{P}(\bar{O})}(1)$ and 
$\eta$ coincides with the pullback of the natural contact structure 
on $\mathbf{P}(O)$. Then $$f: Z \to \mathbf{P}(\bar{O}) \subset  
\mathbf{P}(\g^*)$$ is a morphism determined by the linear system 
corresponding to $V \subset H^0(Z, L)$.} 
\vspace{0.2cm}

{\em Proof}. Let $Z_0 := f^{-1}(\mathbf{P}(O))$ and put $f_0 := f\vert_{Z_0}$. 
There is a commutative diagram 
\begin{equation} 
\begin{CD} 
H^0(Z, \Theta_Z) @>>> H^0(Z, L) \\ 
@V{\iota}VV @V{\cong}VV \\ 
H^0(Z_0, f_0^*\Theta_{\mathbf{P}(O)}) @>>> H^0(Z_0, f_0^*O_{\mathbf{P}(O)}(1)) \\
@A{({f_0}^*)_{\Theta}}AA @A{({f_0}^*)_{O(1)}}AA \\
H^0(\mathbf{P}(O), \Theta_{\mathbf{P}(O)}) @>>> H^0(\mathbf{P}(O), O_{\mathbf{P}(O)}(1))      
\end{CD} 
\end{equation}
All vertical maps are injective.  
The coadjoint action of $G$ on $\mathbf{P}(\g^*)$ determines an embedding 
$\g \subset H^0(\mathbf{P}(O), \Theta_{\mathbf{P}(O)})$. 
Let $$0 \to E \to \Theta_{\mathbf{P}(O)} \stackrel{\bar{\eta}}\to 
O_{\mathbf{P}(O)}(1) \to 0$$ be the contact structure on $\mathbf{P}(O)$. 
On the other hand, the embedding $\mathbf{P}(O) \to \mathbf{P}(\g^*)$ 
induces an injection $\g = H^0(\mathbf{P}(\g^*), O_{\mathbf{P}(\g^*)}(1)) \to 
H^0(\mathbf{P}(O), O_{\mathbf{P}(O)}(1))$. One can check that 
$\bar{\eta}(\g) = \g$ (cf. [Na 3], Remark in pp.23, 24).   
Since $f$ is $G$-equivariant, we have $\iota (\g) = ({f_0}^*)_{\Theta}(\g)$  
By the commutative diagram we have $({f_0}^*)_{O(1)}(\g) = V$.  Q.E.D. 
\vspace{0.2cm}

An important point in Lemma is that $V$ is determined independently 
of $f$. In other words, if $f$ is a morphism satisfying the assumption of 
Lemma, then such an $f$ is unique. 

Let us return to the original situation where 
$Z = \mathbf{P}(\Theta_{G/P})$. We put 
$L := O_{\mathbf{P}(\Theta_{G/P})}(1)$. Then 
the maps $\bar{\pi}$ and $\bar{\mu}$ both satisfy the assumption of 
Lemma. Moreover $L$ is $G$-linearized and the following two maps are 
$G$-equivariant: 
$$\bar{\pi}^*: \g = H^0(\mathbf{P}(\g^*), O_{\mathbf{P}(\g^*)}(1)) 
\to H^0(Z, L)$$  
$$\bar{\mu}^*: \g = H^0(\mathbf{P}(\g^*), O_{\mathbf{P}(\g^*)}(1)) 
\to H^0(Z, L)$$  
As proved in Lemma, $\mathrm{Im}(\bar{\pi}^*) = V$ and 
$\mathrm{Im}(\bar{\mu}^*) = V$. We then have a $G$-equivariant linear 
automorphism  $$\varphi: \g \stackrel{\bar{\pi}^*}\to V \stackrel{(\bar{\mu}^*)^{-1}}\to 
\g.$$
As a consequence we get a commutative diagram  
\begin{equation} 
\begin{CD} 
Z @>{id}>> Z \\ 
@V{\bar{\pi}}VV @V{\bar{\mu}}VV \\ 
\mathbf{P}(\g^*) @>{P(\varphi)}>> \mathbf{P}(\g^*)     
\end{CD} 
\end{equation}
     
Therefore $O$ and 
$O'$ are the same orbit and $\bar{\pi}$ can be regarded as the 
projectivized moment map.  Since the projectivized moment map is 
birational onto its image, the moment map itself is also birational onto 
its image. Thus the moment map (or Springer map) $\mu: T^*(G/P) \to \bar{O}$ 
gives a crepant resolution of $\bar{O}$. 

In order to relate $\mu$ with the original crepant resolution $\pi: Y \to \bar{O}$, 
we need the following lemma. 

{\bf Lemma}. {\em Let $Z$ be the projectivized cotangent bundle 
$\mathbf{P}(\Theta_{G/P})$ of a rational homogeneous space such that 
the projectivized Springer map $\bar{\mu}$ is a birational morphism. 
Denote by $p: Z \to G/P$ the projection map and put 
$L := O_{\mathbf{P}(\Theta_{G/P})}(1)$. Then the nef cone 
$\overline{\mathrm{Amp}}(Z)$ of $Z$ is the closed convex cone generated 
by $[L]$ and  $p^*\overline{\mathrm{Amp}}(G/P)$   
except when $G/P = \mathbf{P}^n$. } 
\vspace{0.2cm}

{\em Proof}. 
Lemma clearly holds true when $b_2(G/P) =1$ and $L$ is not ample. 
We assume that $b_2(G/P) > 1$.  
Note that $\overline{\mathrm{Amp}}(G/P)$ is a simplicial 
polyhedral cone, where each codimension-one face $\mathcal{F}$ corresponds to a morphism $G/P \to 
G/\bar{P}$ with some parabolic subgroup $\bar{P}$ containing $P$.  
Let $\bar{p}: Z \stackrel{p}\to G/P \to G/\bar{P}$ be the 
composed map. 
Then the projectivized Springer map $\bar{\mu}: Z \to \mathbf{P}(\bar{O})$ together 
with $\bar{p}$ gives a morphism $\phi: Z \to G/\bar{P} \times \mathbf{P}(\bar{O})$. 
We prove that $\phi$ actually contracts some curve to a point if $L$ is not 
ample. Let $C$ be a smooth rational curve on $G/P$ contained in a fibre $F$ of the map $G/P \to G/\bar{P}$. 
We have natural surjections $\Theta_{G/P}\vert_C \to N_{C/(G/P)}$ and 
$N_{C/(G/P)} \to N_{F/(G/P)}\vert_C$. 
As $N_{F/(G/P)}\vert_C \cong O_C^{\oplus m}$ 
some $m > 0$, there is a surjection $\Theta_{G/P}\vert_C \to O_C^{\oplus m}$. 
In particular, $\Theta_{G/P}\vert_C$ is not an ample vector bundle. 
This means that $(L.D) = 0$ for some curve $D \subset Z$ with $p(D) = C$. 
Then $D$ is contracted to a point by $\phi$. 
Let $\tilde{\mathcal F}$ be a convex cone generated by $[L]$ and $p^*\mathcal{F}$. 
The observation above shows that $\tilde{\mathcal F}$ is a codimension-one face  
of $\overline{\mathrm{Amp}}(Z)$.   
Finally note that $L$ is ample if and only if  $G/P = \mathbf{P}^n$ by Mori [Mo]. Q.E.D. 
\vspace{0.2cm} 

By the construction of  $Z = Y - \pi^{-1}(0)/\mathbf{C}^*$, some $\pi$-ample 
line bundle $\mathcal{L}$ on $Y$ descends to a $\bar{\pi}$-ample line bundle 
$\bar{\mathcal{L}}$ on $Z$. Assume that $M = G/P$ is not a projective space. 
Then by the lemma above we may assume that    
$\bar{\mathcal{L}} = p^*F$ for some ample line bundle $F$ on $M = G/P$. 
We note that $\pi$ and $\mu$ coincide over $\bar{O} - \{0\}$. 
In fact, since $\bar{\pi}^*O_{\mathbf{P}(\bar{O})}(-1) = O_{\mathbf{P}(\Theta_M)}(-1)$, 
we have $T^*M - \{0-section\} = Z \times_{\mathbf{P}(\bar{O})} (\bar{O} - \{0\})$. 
By the commutative diagram 
\begin{equation} 
\begin{CD} 
Y - \pi^{-1}(0) @>>> Z \\ 
@V{\pi}VV @V{\bar{\pi}}VV \\ 
\bar{O} - \{0\} @>>> \mathbf{P}(\bar{O})     
\end{CD} 
\end{equation}
there is a morphism $Y - \pi^{-1}(0) \to T^*M - \{0-section\}$ over $\bar{O} - \{0\}$. 
Since both $Y - \pi^{-1}(0)$ and $T^*M - \{0-section\}$ are crepant resolutions of 
$\bar{O} - \{0\}$, the morphism is an isomorphism. 
As a consequence we have the commutative diagram 

\begin{equation} 
\begin{CD} 
Y @<{j}<< T^*(G/P) - (0-section) @>{j'}>> T^*(G/P)\\ 
@V{\pi}VV @VVV @V{\mu}VV \\ 
\bar{O} @<<< \bar{O}- \{0\}  @>>> \bar{O}      
\end{CD} 
\end{equation}

Let $q: T^*(G/P) - \{0-section\} \to Z (= \mathbf{P}(\Theta_{G/P}))$ 
be the quotient map. By the definition $j_*q^*\bar{\mathcal{L}} = \mathcal{L}$. 
On the other hand, we have ${j'}_*q^*\bar{\mathcal{L}} = \tilde{p}^*F$, where 
$\tilde{p}: T^*(G/P) \to G/P$ is the projection map. Since $F$ is an ample 
line bundle on $G/P$, $\tilde{p}^*F$ is a $\mu$-ample line bundle. 
The birational map $Y ---> T^*(G/P)$ (over $\bar{O}$) is an isomorphism 
in codimension one and $\tilde{p}^*F$ is the proper transform of 
$\mathcal{L}$ by this map. As each line bundle is $\pi$-ample or $\mu$-ample, 
we see that this birational map is actually an isomorphism.     

We next consider the case when $G/P = \mathbf{P}^n$. 
In this case the contact projective manifold $Z$ has two different 
projectivized cotangent bundle structures over $\mathbf{P}^n$. 
In fact $Z$ is a hypersurface 
of $\mathbf{P}^n \times \mathbf{P}^n$ of type $(1,1)$ and two different 
projections $p_i: \mathbf{P}^n \times \mathbf{P}^n \to \mathbf{P}^n$, 
$i = 1,2$ induce mutually different identifications 
$Z \cong \mathbf{P}(\Theta_{\mathbf{P}^n})$.  
In each choice the corresponding parabolic subgroup $P$ is 
not conjugate to one another. This phenomenon prevents us from recovering 
$Y$ from $Z$. But in this case we can easily check that $\bar{O}$ has 
exactly two different crepant resolution and both of them are Springer 
resolutions.  

Finally we notice that one always has $b_2(Z) > 1$ with only one exceptional 
case when $O$ is the minimal nilpotent orbit of $sl_2$. 
Let $\nu: \tilde{O} \to \bar{O}$ be the normalization. Then $\nu^{-1}(0)$ 
consists of one point $0'$. In fact, the central fibre of 
the Jacobson-Morozov resolution of $\bar{O}$ is connected and the resolution 
factors through $\tilde{O}$; hence $\nu^{-1}(0)$ is one point. 
The $\mathbf{C}^*$-action 
on $\bar{O}$ naturally extends to that on $\tilde{O}$. 
Suppose that $\mathbf{P}(\tilde{O})$ is not smooth. Then the 
crepant resolution $Z \to \mathbf{P}(\tilde{O})$ has exceptional locus; 
in particular, $b_2(Z) > 1$. 
Suppose to the contrary that $\mathbf{P}(\tilde{O})$ is smooth. 
Assume that $\tilde{O}$ is smooth. Then $T_{0'}\tilde{O}$ admits a symplectic 
2-form $\omega$ of weight $1$, which is absurd.  
Thus $\tilde{O}$ has an isolated singularity at $0'$ and 
$Z = \mathbf{P}(\tilde{O})$. 
When $O$ is the minimal nilpotent orbit of $sl_2$, we have $Z = \mathbf{P}(\tilde{O}) 
= \mathbf{P}^1$ and $b_2(Z) = 1$. Otherwise $\dim O \geq 4$. Then the exceptional locus 
of the crepant resolution $\pi: Y \to \tilde{O}$ has codimension $\geq 2$. 
This means that $\tilde{O}$ is not $\mathbf{Q}$-factorial. Note that $Y$ has a 
$\mathbf{C}^*$-action. We take a $\mathbf{C}^*$-linearized $\pi$-ample 
line bundle $L$ on $Y$. We put $M := \Gamma (Y, L)$. Then $M$ has a 
$\mathbf{C}^*$-action. Let $A$ be the coordinate ring of the affine variety 
$\tilde{O}$. Then $A$ is a graded algebra and $M$ is a graded $A$-module. 
Let us consider the coherent sheaf $E := \tilde{M}$ on $\mathbf{P}(\tilde{O}) 
= \mathrm{Proj}(A)$. Then the double dual $E^*$ of $E$ is an invertible sheaf on 
$\mathbf{P}(\tilde{O})$. On the other hand, let $H$ be the pullback of 
$O_{\mathbf{P}(\bar{O})}(1)$ by the map $\mathbf{P}(\tilde{O}) \to 
\mathbf{P}(\bar{O})$. Then $[E^*]$ and $[H]$ are linearly independent in 
$\mathrm{Pic}(\mathbf{P}(\tilde{O}) \otimes \mathbf{Q}$.   
\vspace{0.3cm}

\begin{center}
{\bf References} 
\end{center}

[BCHM] Birkar, C., Cascini, P., Hacon, C., McKernan, J.: 
Existence of minimal models for varieties of log general type, 
J. Amer. Math. Soc. {\bf 23} (2010), 405-468 

[Be] Beauville, A.: Symplectic singularities, Invent. Math. {\bf 139} (2000), 
541-549

[Be 2] Beauville, A.: Fano contact manifolds and nilpotent orbits, 
Comment. Math. Helv. {\bf 73} (1998), 566-583 

[CG] Chriss, N., Ginzburg, V.: Representation theory and complex geometry, 
(1997) Birkhauser 

[CM] Collingwood, D., McGovern, W.: 
Nilpotent orbits in semi-simple Lie algebras, van 
Nostrand Reinhold, Math. Series, 1993

[DPS] Demailly, J.P. , Peternell, T., Schneider, M..: 
Compact complex manifolds with numerically effective 
tangent bundles, J. Algebraic Geom. {\bf 3}. (1994), 295-345  

[F] Fu, B.: Symplectic resolutions for nilpotent orbits, Invent. Math. 
{\bf 151} (2003), 167-186  

[Fu] Fujiki, A.: On automorphism groups of compact Kaehler manifolds, Invent. Math. 
{\bf 44} (1978), 225-258 

[F-H] Fulton, W., Harris, J.: Representation theory, GTM 129 (1991), Springer-Verlag  

[GKK] Greb, D., Kebekus, S., Kov\'{a}cs, S.: 
Extension theorems for differential forms and Bogololov-Sommese vanishing 
on log canonical varieties, Compositio Math. {\bf 146} (2010), 193-219   


[Ka] Kaledin, D.: Symplectic varieties from the Poisson point of view, 
J. Reine Angew. Math. {\bf 600} (2006), 135-160 

[KPSW] Kebekus, S., Peternell, T., Sommese, A., Wisniewski, J.: 
Projective contact manifolds, Invent. Math. {\bf 142}. (2000), 1-15 

[K-P] Kraft, H., Procesi, C.: On the geometry of conjugacy classes 
in classical groups, Comment. Math. Helv. {\bf 57} (1982) 539-602

[LeB] LeBrun, C.: Fano manifolds, contact structures, and quartanionic 
geometry, Intern. J. Math. {\bf 6} No.3 (1995) 419-437 

[LNSV] Lehn, M., Namikawa,Y., Sorger, C., van Straten, D.: 
On symplectic hypersurfaces, preprint, December, 2011 

[Mo] Mori, S.: Projective manifolds with ample tangent bundles, Ann. Math. 
{\bf 110} (1979), 593-606

[Mu] Mumford, D., Forgaty, G., Kirwan, F.: 
Geometric invariant theory, 3-rd enlarged version, 
Ergebnisse der Mathematik und ihrer Grenzgebiete 34, Springer verlag 

[Na 1] Namikawa, Y.: Flops and Poisson deformations of symplectic varieties, 
Publ. Res. Inst. Math. Sci. {\bf 44} (2008), 259-314 

[Na 2] Namikawa, Y.: Poisson deformations of affine symplectic varieties, 
Duke Math. J. {\bf 156} (2011) 51-85 

[Na 3] Namikawa, Y.: Equivalence of symplectic singularities, arXiv: 1102.0865 

[Na 4] Namikawa, Y.: Poisson deformations of affine symplectic varieties II, 
Kyoto J. Math. {\bf 50} No.4 (2010), 727-752 

[Pro] Procesi, C.: Lie groups: an approach through invariants and representations, 
(2007), UTM  Springer-Verlag  

[Sl] Slodowy, P.: Simple singularities and simple algebraic groups, Lecture Notes 
in Mathematics, {\bf 815} (Springer, New York, 1980)    

[Spr] Springer, T.: Linear algebraic groups, 2-nd edition (1998),  Birkh\"{a}user

\vspace{0.4cm}

Department of Mathematics, Faculty of Science, Kyoto University 

e-mail address: namikawa@math.kyoto-u.ac.jp

\end{document}